\documentclass[11pt]{elsarticle}
\usepackage{amsmath}
\usepackage{amssymb}
\usepackage{amsthm}
\usepackage{mathrsfs}
\usepackage[nodots]{numcompress}

\theoremstyle{plain} 
\newtheorem{cor}{Corollary}

\newtheorem{pro}{Proposition}
\newtheorem{thm}{Theorem}

\theoremstyle{definition}

\theoremstyle{remark}

\newcommand{\set}[1]{\left\{#1\right\}}
\newcommand{\pd}{\,\partial}
\newcommand{\comment}[1]{}

\begin{document}

\begin{frontmatter}



\title{Coefficient characterization  of linear  differential equations with maximal symmetries}

\author{J.C. Ndogmo\corref{cor1}}
\ead{jean-claude.ndogmo@wits.ac.za}


\address{School of Mathematics,
University of the
Witwatersrand,
Private Bag 3, Wits 2050,
South Africa}

\begin{abstract}
A characterization of the general linear equation in standard form admitting a maximal symmetry algebra is obtained in terms of a simple set of conditions relating the coefficients of the equation. As a consequence, it is shown that in its general form such an equation can be expressed in terms of only two arbitrary functions, and its connection with the Laguerre-Forsyth form is clarified. The characterizing conditions are also used to derive an infinite family of semi-invariants, each corresponding to an arbitrary order of the linear equation. Finally a simplifying ansatz is established, which allows an easier determination of the infinitesimal generators of the induced pseudo group of equivalence transformations, for all the three most general canonical forms of the equation.
\end{abstract}


\end{frontmatter}

\section{Introduction}
\label{s:intro}

By a result of Lie \cite{liemax}, a linear ordinary differential equation ({\sc ode}) of a general order $n$
is known to have a symmetry algebra of maximal dimension $d_n$ if it is reducible by a point transformation to the
equation $ y^{(n)}=0$, which will henceforth be referred to as the canonical form of the linear equation. In a much recent
paper Krause and Michel \cite{KM1} proved the converse of this result and also showed that a linear equation is iterative
if and only if its symmetry algebra has the maximal dimension $d_n.$ (By the cited result of Lie \cite{liemax}, $d_n= n+4$ for $n\geq 3$). Characterizing linear equations having a symmetry algebra
of maximal dimension is therefore the same as characterizing linear equations that are reducible by a point transformation
to the canonical form. The latter characterization for the third-order equation $ y^{(3)} + c_2\, y'' + c_1\,y' + c_0\, y=0 $
is due to Lie \cite{lieGc} and Laguerre \cite{Lag3} who showed independently that this equation is reducible to the canonical
form if and only if its coefficients satisfy the equation

\begin{equation}\label{eq:1.1}
 54 c_0 - 18 c_1 c_2 + 4 c_2^3 - 27 c_1' + 18 c_2 c_2' + 9 c_2''=0.
\end{equation}

     This characterization also clearly applies to all nonlinear  {\sc ode}s which are linearizable  by point transformations \cite{faz-leach, ibra-magri}, as the latter transformations do not alter the dimension of the symmetry algebra.\par
     In this paper, we extend this characterization to equations of higher orders. It turns out that for each equation
of order $n$ there will be $n-2$ characterizing equations, and the limitation of our presentation of the characterizing equations
only up to the order five is simpy due to their very large size. However, we give a description of the method for deriving this
characterization for equations of any order. The derivation of these characterizing equations is also based on the canonical
 normal form of linear equations admitting a maximal symmetry algebra that was obtained in \cite{faz-leach} from a symmetry approach,
and in \cite{jcfaz} from an iterative approach. These characterizing equations therefore also represent a generalization of the results
of \cite{faz-leach} and \cite{jcfaz}. We then deduce that the most general form of a linear equation admitting a maximal symmetry algebra can be expressed in standard form in terms of only two arbitrary functions. We also deduce that the Laguerre-Forsyth form  of a linear equation reduces to the canonical form if and only if the equation has maximal symmetries.

     Although we do not give the characterizing equations for each linear equation of order $n$, we note however that among the
$n-2$  characterizing equations exactly one of them represents a semi-invariant of the equation, that is a function of the coefficients
of the equation whose expression does not change when the dependent variable is transformed. We obtain an expression for these semi-invariants for equations of all orders and describe some of their properties.

Finally, using some simplifying assumptions and the method of \cite{jcftc}, we give expressions for both the symmetry generator $X_n$ of $G_{\!S}$  and $X_n^0$ of the induced pseudo group of transformations $G_c,$  and for all three  most general canonical forms of linear equations of  a general order $n.$ Here,  $G_{\!S}$ denotes the symmetry group of the general linear equation in which the arbitrary functions are considered as additional dependent variables.

\section{Coefficient characterization}
    A method based on a symmetry approach has been proposed in \cite{faz-leach} for characterizing the coefficients of linear ordinary differential equations ({\sc ode}s) that admit a maximal symmetry algebra, but only for equations in reduced normal form (in which the term of second highest order vanishes). In  a more recent paper \cite{jcfaz} a similar characterization based on an iterative approach was proposed, in which according to a result of Krause and Michel \cite{KM1} a linear equation admitting a maximal symmetry is simply viewed as an iterative equation. By iterative equation, we mean an equation of the form
\begin{subequations}\label{eq:2.1}
\begin{align}
\Psi ^n[y]&=0, \quad y=y(x), \quad n \geq 1 \\
\intertext{where}
\vspace{-5mm} \Psi^1 [y]&= r y' + s y, \quad \Psi^n [y]=\Psi^{n-1}\left[ \Psi[y]\right],
\end{align}
\end{subequations}
 and where $r=r(x)$ and $s=s(x)$ are the parameters of the source equation $\Psi^1[y]=0$. This characterization shows that in its reduced normal form, a general linear equation depends solely on one arbitrary function $a = a(x)$. For equations of orders three to five, the corresponding equations are given as follows:
\begin{subequations}\label{eq:2.2}
\begin{align}
&y^{(3)} + a y' + \frac{a'}{2} y =0 \label{eq:2.2a}\\
&y^{(4)} + a y'' + a' y' + \left(\frac{3}{10} a'' + \frac{9}{100} a^2\right)y =0 \label{eq:2.2b}\\
&y^{(5)} + a y^{(3)} + \frac{3}{2} a' y'' + \left(  \frac{9}{10} a'' + \frac{16}{100} a^2 \right) y' + \left( \frac{1}{5} a^{(3)} + \frac{16}{100} a a'  \right) y=0.
\end{align}
\end{subequations}

    However, as a linear equation need not occur in its reduced normal form, but rather in the most general standard form, it is thus useful to obtain the corresponding characterization for equations in standard form. We let the general linear equation be given in standard form as
\begin{equation} \label{eq:2.3}
\Delta(x, y_{(n)}; C)\; \equiv \;y^{(n)} + c_{n-1}\, y^{(n-1)} + c_{n-2}\, y^{(n-2)} + \dots + c_0\, y=0.
\end{equation}
where $C= (c_0, \dots, c_{n-1}).$ Suppose that such an equation has a symmetry algebra of maximal dimension and let its corresponding reduced normal form be given by
\begin{equation} \label{eq:2.4}
y^{(n)} + B_{n-2}\, y^{(n-2)} + B_{n-3}\, y^{(n-3)} + \dots + B_0\, y=0,
\end{equation}

 where the $ B_j $ for $ j = 0,\dots , n-2 $ are its coefficients and depend as already noted
above on a single arbitrary function $ a = B_{n-2} $ and its derivatives. Let
\begin{equation} \label{eq:2.5}
y^{(n)} + A_{n-1}\, y^{(n-1)} + A_{n-2}\, y^{(n-2)} + \dots + A_0\, y=0
\end{equation}
be the corresponding standard form of \eqref{eq:2.4}, which may be obtained by a transformation of the form
\begin{equation} \label{eq:2.6}
y \mapsto y e^{- \frac{1}{n} \int_{x_0}^x A_{n-1} dx}.
\end{equation}

Then \eqref{eq:2.3} and \eqref{eq:2.5} must be identical, and in particular the nonzero coefficient $A_{n-1}$ introduced by the transformation \eqref{eq:2.6} satisfies $A_{n-1}~= ~c_{n-1}$, and more generally we have
\begin{equation} \label{eq:2.7}
c_j = A_j,\quad \text{ for $j=0, \dots, n-1$}.
\end{equation}
Note that the coefficients $c_j$ in \eqref{eq:2.3} are mere symbols and we wish to find a relationship among them. Given that in \eqref{eq:2.4} the function $B_{n-2}$ is precisely the arbitrary function $a(x)$ labeling the equation, it can be shown by a recursive procedure,  or even by induction on $n$ that
$$
A_{n-2} = a + \frac{n-1}{2n}\, c_{n-1}^2 + \frac{n-2}{2}\, c_{n-1}'.
$$
Therefore, solving the equation $ c_{n-2} = A_{n-2} $ for $a$ gives
\begin{equation} \label{eq:2.8}
a= c_{n-2} - \left( \frac{n-1}{2n} \, c_{n-1}^2 + \frac{n-2}{2}\, c_{n-1}'\right).
\end{equation}
    Consequently, the characterizing equations for linear equations in standard form with maximal symmetry algebra are given by the remaining $n-2$ equations
\begin{equation} \label{eq:2.9}
c_j = A_j, \quad \text{$j=0, \dots, n-3$},
\end{equation}
in which the function $a$ and its derivatives are substituted with the corresponding expressions given by \eqref{eq:2.8}.
\begin{pro}\label{pro:2.1}
If a linear equation in standard form \eqref{eq:2.3} has maximal symmetry, then in its general form it may be expressed in terms of only two arbitrary functions, namely the functions $c_{n-1}$ and $c_{n-2}$, and their derivatives.
\end{pro}
\begin{proof}
The result readily follows from the fact that the functions $A_j$ in \eqref{eq:2.9} then depend only on $a$ and its derivatives, while \eqref{eq:2.8} shows that the function $a$ depends precisely on $c_{n-1}$, $c_{n-2}$, and their derivatives.
\end{proof}
\begin{cor}
A linear equation in standard form \eqref{eq:2.3} with $c_{n-1}=c_{n-2}=~0$ has maximal symmetry algebra if and only if $c_j=0$ for all $j$. In other words a linear equation has maximal symmetry algebra if and only if its Laguerre-Forsyth form corresponds to the canonical equation $y^{(n)}=0$.
\end{cor}
\begin{proof}
After all a Laguerre transformation is also a point transformation although it cannot always be explicitly constructed for a given equation. Since equations equivalent under point transformation have similar Lie algebras, it readily follows that if the Laguerre-Forsyth form of an equation is $y^{(n)}=0$, then the equation has maximal symmetry algebra. The converse of the corollary is a direct application of proposition \ref{pro:2.1}, and the fact that in \eqref{eq:2.9} the $c_j$ turn out to be polynomial functions with no constant terms of $c_{n-1}$, $c_{n-2}$, and their derivatives.
\end{proof}

As an immediate consequence of the corollary, linear equations such as $y^{(3)} + f(x)y=0$ or $y^{(4)} +f(x)y'=0$ have maximal symmetry algebras if and only if the function $f(x)$ vanishes identically.
We now make use of \eqref{eq:2.9} and \eqref{eq:2.8} to explicitly derive the characterizing equations for maximal symmetry algebras for equations of orders three to five.\par

     For $n=3$, it is readily found that in \eqref{eq:2.5} we have
\begin{equation} \label{eq:2.10}
A_0 = \frac{1}{54} \left( 18 a c_2 + 2 c_2^3 + 27 a' + 18 c_2 c_2' + 18 c_2''\right),
\end{equation}
while the corresponding expression for $a$ in \eqref{eq:2.8} reduces to
\begin{equation} \label{eq:2.11}
a= c_1 - \left( \frac{c_2^2}{3} + \frac{c_2'}{2}\right).
\end{equation}
Applying \eqref{eq:2.11} into \eqref{eq:2.10} and substituting the resulting expression for $A_0$ into \eqref{eq:2.9} gives exactly the already cited equation \eqref{eq:1.1} found by Lie \cite{lieGc} and Laguerre \cite{Lag3} and given by
$$54c_0 - 18c_1c_2 + 4 c_2^3 - 27c_1' + 18c_2c_2' + 9c_2''=0.$$
   The most general form of a linear third-order equation admitting a maximal symmetry algebra can thus be expressed in terms of only two arbitrary functions $c_1(x)$ and $c_2(x)$ in the form of
\begin{equation} \label{eq:2.12}
y^{(3)} + c_2\, y'' + c_1\, y' + \frac{1}{54} \left(  18 c_1 c_2 - 4 c_2^3 + 27 c_1' -18 c_2 c_2' - 9 c_2''\right)y=0.
\end{equation}
Equation \eqref{eq:2.12} naturally reduces to \eqref{eq:2.2a} for $c_2=0.$\par
For $n=4,$ we successively get
\begin{subequations}\label{eq:2.13}
\begin{align}
a &= \frac{1}{8} \left(  8 c_2 - 3 c_3^2 - 12 c_3' \right) \label{eq:2.13a} \\
A_1 &= \frac{1}{2} \left(  a c_3 + \frac{c_3^3}{16} + a' + \frac{3}{4} c_3 c_3' + c_3'' \right) \label{eq:2.13b}\\
\begin{split}\label{eq:2.13c}
6400 A_0 &=  576 a^2 + 400 a (c_3^2 + 4 c_3')\\
& + 5 \left[  5 c_3^4 + 120 c_3^2 c_3' + 320 c_3 (a' + c_3'')\right]\\
&+ 80 \left( 15 c_3'^2 + 24 a'' + 20 c_3^{(3)}\right).
\end{split}
\end{align}
\end{subequations}
Substituting \eqref{eq:2.13a} into \eqref{eq:2.13b} and \eqref{eq:2.13c} gives the two equations
\begin{subequations}\label{eq:2.14}
\begin{align}
8 c_1 +  &= 4 c_2 c_3 - c_3^3 + 8 c_2' - 6 c_3c_3' - 4 c_3''  \label{eq:2.14a} \\
\begin{split}\label{eq:2.14b}
1600 c_0  &= 144 c_2^2 - 11 c_3^4 +400c_3c_2' - 288 c_3^2 c_3'- 336 c_3'^2 \\
  &\quad - 8 c_2(c_3^2 + 4 c_3') +480 c_2'' - 560 c_3 c_3'' - 320 c_3^{(3)}
\end{split}
\end{align}
\end{subequations}
which represent the characterizing equations for maximal symmetry algebra for equations of order $4$. Note that conversely any linear fourth order equation whose coefficients satisfy \eqref{eq:2.14} must be iterative, which is why conditions such as \eqref{eq:2.14} are termed characterizing equations. Indeed,
if the coefficients of a fourth order equation of the form \eqref{eq:2.3} satisfy \eqref{eq:2.4}, then its reduced normal form has, after the substitution of the expressions for $c_0$ and $c_1$ given by \eqref{eq:2.4} in terms of $c_2, c_3,$ and their derivatives, the form
\begin{subequations}\label{eq:nor4std4}
\begin{align}
w^{(4)} & + Q_2 w'' + Q_1 w' + Q_0 w  =0 \label{eq:nor4std4E} \\
\intertext{where}
Q_2 &= c_2 - \frac{3}{8}(c_3^2 + 4 c_3')\\
Q_1 &= c_2'- \frac{3}{4}(c_3 c_3'+ 2 c_3'')\\
 \begin{split}  Q_0 &= \frac{3}{6400}(192 c_2^2 + 27 c_3^4 - 48 c_3'^2 - 144 c_2 (c_3^2+4 c_3'))\\
&\quad + \frac{3}{6400}( 27 c_3^4+ 216 c_3^2 c_3'+ 640 c_2''- 480 c_3 c_3'' -960 c_3'''). \end{split}
\end{align}
\end{subequations}
The coefficients $Q_j$ thus obtained clearly satisfy the conditions
$$Q_1= Q_2'\quad \text{ and }\quad  Q_0= (\frac{3}{10} Q_2'' + \frac{9}{100} Q_2^2)$$
prescribed by \eqref{eq:2.2b} for iterative equations, as required.
\par

For equations of order $n=5$, by proceeding as above for the orders three and four, we obtain the following $n-2=3$ characterizing equations
\begin{subequations}\label{eq:2.16}
\begin{align}
c_2 &=  (30 c_3 c_4 - 8 c_4^3 +75 c_3' - 60 c_4c_4' + 50 c_4'')/50 \label{eq:2.16a} \\[3mm]
\begin{split} \label{eq:2.16b}
1250\, c_1 &=  + 200 c_3^2 - 18 c_4^4 + 750 c_4 c_3' - 580 c_4^2 c_4'- 850 c_4'^2  \\
 &\quad - 10 c_3 (c_4^2 + 5 c_4') + 1125 c_3'' - 1400 c_4 c_4'' - 1000 c_4^{(3)}
\end{split}\\[3mm]
\begin{split}
6250\,c_0 &= 200 c_3^2 c_4 +14 c_4^5 - 25 c_4^2c_3' + 40c_4^3 c_4' \\
&\quad - 125 c_3'c_4' - 750 c_4 c_4'^2 + 1125 c_4 c_3'' - 850 c_4^2 c_4'' \\
&\quad - 2750 c_4'c_4'' +1250 c_3^{(3)} - 2000 c_4 c_4^{(3)} - 1250 c_4^{(4)}\\
&\quad - 10 c_3 (11 c_4^3 +100c_3' - 85 c_4 c_4' - 75 c_4'').
\end{split}
\end{align}
\end{subequations}
\section{Semi-invariants of linear equations}

      The group of equivalence transformations of the general linear equation \eqref{eq:2.3} is given by invertible point transformations of the form
\begin{equation} \label{eq:3.1}
x=f(z),\quad    y=g(z) w(z),
\end{equation}
and they preserve the linearity and the homogeneity of the equation. Let
\begin{equation} \label{eq:3.2}
w^{(n)} + Q_{n-1}\,w^{(n-1)} + Q_{n-2}\,w^{(n-2)} + \dots + Q_0\, w=0
\end{equation}
be the transformed version of \eqref{eq:2.3} under \eqref{eq:3.1}. By a semi-invariant of \eqref{eq:2.3} we shall mean a function $F= F(c_0, c_1,\dots, c_{n-1})$ of the coefficients of the equation which have the same expression for the transformed equation when the dependent variable (alone) changes. It is well known that under \eqref{eq:3.1} the expression of the semi-invariant for the transformed equation is related to that for the original equation \cite{jMF, forsy} by the equality
\begin{equation} \label{eq:3.3}
F(Q_0, Q_1, \dots, Q_{n-1})= \left( \frac{dx}{dz} \right)^\mu F(c_0, c_1, \dots, c_{n-1}),
\end{equation}
where $\mu$ is an integer. In this case we say that the semi-variant $F$ has index $\mu.$  To each expression of the form $d^k c_j/ d x^k,$
let us assign the weight $(n-j) + k$, and we let this weight function be multiplicative so that the product
$c_p c_q$ has weight $(n-p)+ (n-q).$
It is well known that for a given semi-invariant all terms have the same weight and that this weight coincides with the index of the semi-invariant \cite{jMF, forsy} .\par

A closer look at the set of characterizing equations \eqref{eq:2.9} shows that precisely one of them corresponds to a semi-invariant of the equation, namely the relation $c_{n-3}=A_{n-3}$, which gives rise to the semi-invariant $F=A_{n-3} - c_{n-3}$.\par

First of all, using the method of either \cite{jcfaz} or \cite{faz-leach}, it can be proved that the coefficient $B_{n-3}$ in \eqref{eq:2.4} satisfies
$B_{n-3}= \frac{n-2}{2} a'.$ Consequently, using the expression of the function $a$ in \eqref{eq:2.8} it follows by induction on $n$ that the coefficient $A_{n-3}$ in \eqref{eq:2.5} is given by
\begin{equation} \label{3.4}
\begin{split}
A_{n-3} = &\frac{n-2}{n} c_{n-1} c_{n-2} -\frac{(n-1)(n-2)}{3 n^2} c_{n-1}^3 + \frac{n-2}{2} c_{n-2}'\\
          & - \frac{(n-1)(n-2)}{2n} c_{n-1} c_{n-1}'- \frac{(n-1)(n-2)}{12} c_{n-1}'',
\end{split}
\end{equation}
so that the corresponding invariant function $I_n$ has expression
\begin{equation} \label{eq:3.5}
\begin{split}
I_n = &\frac{n-2}{n} c_{n-1} c_{n-2} -\frac{(n-1)(n-2)}{3 n^2} c_{n-1}^3 + \frac{n-2}{2} c_{n-2}'\\
          & - \frac{(n-1)(n-2)}{2n} c_{n-1} c_{n-1}'- \frac{(n-1)(n-2)}{12} c_{n-1}'' - c_{n-3}.
\end{split}
\end{equation}
The fact that the function $I_n=I_n (c_0, c_1,\dots, c_{n-1})$ in \eqref{eq:3.5} is a semi-invariant can readily be verified. First each term in this expression has weight three, and we readily see that
\begin{equation*}
I_n (Q_0, Q_1, \dots, Q_{n-1})= f'(z)^3 I_n (c_0, c_1, \dots, c_{n-1}),
\end{equation*}
which proves the assertion.\par

    Although the invariant functions $I_n$ in \eqref{eq:3.5} are originally defined only for
$n \geq 3,$   their expression shows that they vanish identically for $n=1$ or $n=2$, by letting $c_j=0$ for $j<0.$
This vanishing can be interpreted by the fact that all first order and all second order linear equations are all equivalent through a point transformation to the equations $y'=0$ and $y''=0$, respectively, and therefore they do not have nontrivial invariant functions.\par

    On the other hand it should be noted that the other equations in the characterizing system \eqref{eq:2.9} do not give rise to invariant functions except for the value $j=n-3$ in that system of equations. Indeed, denote collectively by $C$ and $Q$ the coefficients in equations \eqref{eq:2.3} and \eqref{eq:3.2}, respectively, and for $n=4$ denote by $J(C)=1600 (c_0 - A_0)$ the normalized function obtained from 2.9 with $j=0$. Then it can be seen that although each term in the expression of $J(C)$ has weight four, we have
\begin{equation*}
J(Q)= f'(z)^4 J(C) - 200 \frac{h'(z)}{h(z)} f'(z)^3 I_4 (C),
\end{equation*}
clearly showing that the function $J$ is not a semi-invariant.

\section{Infinitesimal generators of the induced group action}
\label{sec:infigen}
         The equivalence group $G$ in \eqref{eq:3.1} of the general linear equation \eqref{eq:2.3} induces another Lie pseudo group $G_c$ acting on the coefficients of  \eqref{eq:2.3} \cite{lieGc}. For linear equations with maximal symmetries, their most general form depends as already noted on only two arbitrary functions, instead of $n$. For instance, the most general form of linear equations of order four admitting a maximal symmetry algebra is given on account of \eqref{eq:2.14} by
\begin{equation} \label{eq:4.1}
\begin{split}
y^{(4)} & + c_3 y^{(3)} + c_2 y'' + \frac{1}{8} \Big[ 4 c_2 c_3 - c_3^3 + 8 c_2' - 6 c_3 c_3'- 4c_3''\Big] y'\\
&+ \frac{1}{1600} \Big[ 144 c_2^2 -11 c_3^4 + 400 c_3 c_2' -288 c_3^2c_3'- 336c_3'^2 \\
& \qquad - 8c_2 \big(c_3^2 + 4 c_3'\big)+ 480 c_2'' - 560 c_3 c_3'' -320 c_3^{(3)}\Big ] y=0
\end{split}
\end{equation}
and it is expressible solely in terms of the coefficients $c_{n-1}$ and $c_{n-2}$, here $c_3$ and $c_2$.\par

        Although Eq. \eqref{eq:4.1} is a very special case of the general Eq.  \eqref{eq:2.3}, its equivalence group is the same group $G$ in \eqref{eq:3.1} because equivalent equations have similar symmetry groups. Consequently the infinitesimal generators $X^0$ of the group $G_c$ for \eqref{eq:2.3} will also be valid for equations with maximal symmetries. In particular to obtain the specific infinitesimal generators for equations with maximal symmetries expressed only in terms of the two arbitrary functions, it will be sufficient to substitute the characterizing equations \eqref{eq:2.9} into the expression for $X^0$.\par

       A method for finding the infinitesimal generator $X^0$ has been proposed in \cite{jcftc}. If we denote by
\begin{equation} \label{eq:4.2}
X= \xi \pd_x + \eta \pd_y + \phi_{n-1} \pd_{c_{n-1}} + \dots + \phi_0 \pd_{c_0}
\end{equation}
the infinitesimal  generator of \eqref{eq:2.3} in which the coefficients
$$C= (c_0, c_1,\dots, c_{n-1})$$
are also considered as dependent variables, then the method of \cite{jcftc} consists of finding a set of minimum conditions for which the projection $V= \xi\pd_x + \eta\pd_y$ of $X$ on the $(x,y)$-space reduces to the infinitesimal generator $V^0=\set{\xi^0, \eta^0}$ of the equivalence group $G$. This set of minimal conditions imposed to $\phi= (\phi_0, \phi_1,\dots, \phi_{n-1})$ yields a function $\phi^0= (\phi_0^0, \phi_1^0, \dots, \phi_{n-1}^0)$ so that the expression for $X^0$ takes the form
\begin{equation} \label{eq:4.3}
X^0= \xi^0 \pd_x + \phi_{n-1}^0 \pd_{c_{n-1}} + \dots + \phi_0^0 \pd_{c_0}.
\end{equation}
   In practice, the determination of the symmetry generator $X$ for the general linear equation \eqref{eq:2.3} is computationally exhaustive, and a popular Lie symmetry software such as MathLie (See \cite{baumann}) computes $X$ only for $n \leq 4$ due to computer memory problems (on an Intel Core2 Quad CPU machine) while another well-known similar Lie symmetry software such as SYM \cite{sym} does not compute symmetries such as $X$ that involve several dependent variables for a single independent variable.\par

      We therefore need an efficient simplifying ansatz for the manual computation of $X^0$ at orders higher than the fourth. For this, we note that as the full symmetry group of \eqref{eq:2.3} with $C$ considered also as dependent variable should leave the equation invariant, the transformation of the dependent and the independent variables should preserve the form of the equation, except for the introduction of a constant term independent of $y$ which should be offset by the subsequent transformations of the coefficient $C$. This means that in \eqref{eq:4.2}, we must have
\begin{equation} \label{eq:4.4}
\xi= f(x), \qquad \eta= g(x) y + h(x).
\end{equation}
     A verification of \eqref{eq:4.4} is possible by direct calculation for equations of order not higher than the fourth using the MathLie software, while for orders higher than four, the validity of the generators $X$ and $X^0$ found can be tested through the satisfaction of the infinitesimal condition of invariance applied to the general linear equation \eqref{eq:2.3}, and to the semi-invariants $I_n$ found in \eqref{eq:3.5}, respectively. Recall that the infinitesimal criterion of invariance for the infinitesimal generator $X$ of \eqref{eq:2.3} is given by
\begin{equation}\label{eq:infi-cond}
     X^{[n]} \left[  \Delta(x, y_{(n)}; C)\right] =0, \quad \text{whenever $ \Delta(x, y_{(n)}; C) =0,$}
     \end{equation}
where  $X^{[n]}$ represents the $n$-th prolongation of $X.$  Regarding the verification of the infinitesimal condition of invariance for semi-invariants, we note that if for some group element $\alpha \in G_c$ we set $Q= \alpha \cdot C,$ then every semi-invariant of $G_c$ satisfies $F(\alpha \cdot C)= \mathbf{w}(\alpha) \cdot F(C)$  for some weight function $\mathbf{w},$ and $X^0$ is an infinitesimal generator of $G_c$ if and only if
$$
X^0 \cdot F= - d \mathbf{w}(e) F,
$$
for all such functions $F,$ where $\mathbf{w}(e)$ is the differential of $\mathbf{w}$ at the identity element $e$ of $G_c.$ In the actual case of \eqref{eq:2.3} and $G_c$ (which is the same as G except that it acts on the space of coefficients), for $\alpha \equiv (f, g)$ specified in \eqref{eq:3.1} we have $\mathbf{w}(\alpha)= f'(z)^3,$ and for each generator $X^0 \equiv X^0 (n)$ found, it is readily verified that
\begin{equation} \label{eq:4.6Z}
X^0 \cdot I_n = - 3 f'(x) I_n,
\end{equation}
as required.\par

To our knowledge the infinitesimal generators $X^0$ of the induced pseudo group $G_c$ has been computed only for third order equations, or for the normal or the Laguerre-Forsyth  forms of equations of low orders not exceeding five \cite{ibranot, jcpla, jMF}. This is due in part as already mentioned to the intensive computational requirements for the calculation of these generators, but also because the more systematic method for finding them proposed in \cite{jcftc} is relatively recent.


We list in the next three theorems the general expressions for the infinitesimal generators $X_n$ of $G_{\!S}$ and $X_n^0$ of $G_{\!c}$ and for the three most general canonical forms of linear equations, where the subscript $n$ denotes the order of the equation.
\begin{thm}\label{th:std}
For the general linear equation of order $n$ in standard form \eqref{eq:2.3}, the infinitesimal generators $X_n$ of $G_{\!S}$ and $X_n^0$ of $G_{\!c}$ have the following expressions, where $f,g$ and $h$ are arbitrary functions of $x$, and $\delta_0^k$ denotes the Kronecker delta.
\begin{enumerate}
\begin{subequations}\label{eq:4.6}
\item[{\rm a)}]
\begin{align}
& X_n = f \partial_x + \left(yg + h \right)\partial _y + \sum_{k=0}^{n-1} \Phi_k^n \partial _{c_k}, \label{eq:4.6a} \\
 \intertext{where}
\begin{split}
&\Phi _k^n = -(n-k)c_k f' + \sum_{j=1}^{n-k} c_{k+j}\left[\binom{k+j}{j+1}f^{(j+1)}- \binom{k+j}{j}g^{(j)} \right] \\
&\quad \qquad + \delta _0^k \left[-c_k \dfrac{h}{y} + \sum_{j=1}^{n-k} c_{k+j}\binom{k+j}{j}\dfrac{h^{(j)}}{y} \right],\quad  \text{for $k=0, \dots, n-1$.} \label{eq:4.6b}
\end{split}
\end{align}
\end{subequations}
\item[{\rm b)}]
\begin{subequations}\label{eq:4.7}
\begin{align}
& X_n^0 = f \partial_x +\sum_{k=0}^{n-1}\Phi _k^n \partial _{c_k},   \label{eq:4.7b}
\intertext{where}
\begin{split}
&\Phi _k^n = -(n-k)c_k f' + \sum_{j=1}^{n-k} \left[-\binom{k+j}{j}g^{(j)}+\binom{k+j}{j+1}f^{(j+1)} \right] c_{k+j}, \\
&\qquad \qquad \text{for $k=0, \dots, n-1$.}
\end{split}
\end{align}
\end{subequations}
\end{enumerate}
\end{thm}
\begin{proof}
We let the generator $X_n$ be in the form
\begin{align}
X_n =& \xi \partial _x + \eta \partial _y +\sum_{k=0}^{n-1} \Phi _k^n \partial _{c_k},
\end{align}
where the functions $\xi$, $\eta$, and $\Phi _k^n$ are to be specified. We know from the ansatz \eqref{eq:4.4} that $\xi =f(x)$ and $\eta = g(x) y + h(x)$ for some arbitrary functions $f, g$ and $h$ of $x$. The prolongation formula for $X_n ^{[n]}$ is well-known \cite{olver86}. Writing down this expression and applying the infinitesimal condition of invariance \eqref{eq:infi-cond} gives the usual determining equations for the coefficients $\xi, \eta$ and $\Phi_k^n$. Although the procedure is a lengthy  one, thanks to the ansatz \eqref{eq:4.4} these determining equations are easily solved and lead to the expressions in \eqref{eq:4.6}.\par
For the second part of the theorem, the result follows by noting that according to the algorithm of \cite{jcftc} already described for finding $X_n^0$, one essentially only need to find the minimum set of conditions which reduce the projection $\{f(x), g(x)y + h(x) \}$ of $X_n$ onto the $(x,y)$-space to the infinitesimal generator of the equivalence group. From the expressions of the equivalence transformations given in \eqref{eq:3.1}, it follows that the required minimal set of condition reduces to $\{h=0\}$. Applying these conditions to \eqref{eq:4.6} and dropping the term in $\partial _y$ gives the required expression \eqref{eq:4.7}.
\end{proof}
\begin{thm}\label{th:nor}
For the general linear equation in reduced normal form, i.e. in the form \eqref{eq:2.3} with $c_{n-1}=0$, the generators $X_n$ of $G_{\!S}$ and $X_n^0$ of $G_{\!c}$ have the following expressions, in terms of the arbitrary functions $f$ and $h$ of $x$.
\begin{enumerate}
\item[{\rm a)}]
\begin{subequations}\label{eq:4.8}
\begin{align}
& X_n = f \partial_x + \left[y\left[\left(\dfrac{n-1}{2}\right)f'+K_1\right] + h \right]\partial _y + \sum_{k-0}^{n-2} \Phi _k^n \partial _{c_k}, \intertext{where}
\begin{split}
&\Phi _k^n = -(n-k) f' c_k  + \sum_{j=1}^{n-k} c_{k+j}\left[\binom{k+j}{j+1}- \binom{k+j}{j}\dfrac{n-1}{2} \right]f^{(j+1)} \\
&\qquad \quad + \delta _0^k \left[-c_k \dfrac{h}{y} + \sum_{j=1}^{n-k} c_{k+j}\binom{k+j}{j}\dfrac{h^{(j)}}{y} \right],\quad  \text{for $k=0, \dots, n-2$.}
\end{split}
\end{align}
\end{subequations}
\item[{\rm b)}]
\begin{subequations}\label{eq:4.9}
\begin{align}
& X_n^0 = f \partial_x +\sum_{k=0}^{n-2}\Phi _k^n \partial _{c_k}, \intertext{where}
\begin{split}
&\Phi _k^n = -(n-k)c_k f' + \sum_{j=1}^{n-k} a_{k+j}\left[\binom{k+j}{j+1}-\binom{k+j}{j}\dfrac{n-1}{2}\right]f^{(j+1)} , \\
&\qquad \quad \text{for $k=0, \dots, n-2$.}
\end{split}
\end{align}
\end{subequations}
\end{enumerate}
\end{thm}
\begin{proof}
The expressions for $X_n$ and $X_n^0$ are to be sought in the form \eqref{eq:4.6} and \eqref{eq:4.7}, respectively, as the normal form of \eqref{eq:2.3} is a special case of that equation. The main difference is that the equivalence transformations for the normal form are no longer given by \eqref{eq:3.1} but by the much restricted version
\begin{align}\label{eq:4.10}
x=& T(z), \quad y=\lambda \left[ T'(z)\right]^{\frac{n-1}{2}} w(z)
\end{align}
where $T$ is an arbitrary function and $\lambda$ an arbitrary constant. This has infinitesimal generator
\begin{align}\label{eq:4.10'}
V=& f(x) \partial _x + y\left( \dfrac{n-1}{2}f'(x) + k_1\right) \partial _y,
\end{align}
where $f$ is an arbitrary function and $k_1$ an arbitrary constant. Since the functions $f$ and $g$ in \eqref{eq:4.6} and \eqref{eq:4.7} are precisely the parameters of the infinitesimal generator of the equivalence group, to obtain \eqref{eq:4.8} and \eqref{eq:4.9}, we only need to replace $g$ in the latter expressions by the substitution $g=\frac{n-1}{2}f' + k_1$ and to drop the term in $c_{n-1}$. This yields \eqref{eq:4.8} and \eqref{eq:4.9}.
\end{proof}
The Laguerre-Forsyth form of the general linear equation is the equation of the form \eqref{eq:2.3} in which the coefficients $c_{n-1}$ and $c_{n-2}$ of terms of second and third highest orders have vanished. In principle, such a transformation can be realized by means of the change of variables of the form
\begin{subequations}\label{eq:4.11}
\begin{align}
\set{z, x} &= \frac{12}{n(n-1)(n+1)} c_{n-2}, \quad y= w  \exp\left(- \frac{1}{n} \int_{z_0}^z c_{n-1} dx \right),\\
\intertext{ where }
\set{z, x}&= \big[  z' z^{(3)} - (3/2) z''^2 \big] z'^{\,-2}
\end{align}
\end{subequations}
is the Schwarzian derivative, and $z'= dz/dx.$ The Laguerre-Forsyth form of \eqref{eq:2.3} is therefore of an implicit nature in the sense that \eqref{eq:4.11} can not always be solved explicitly for $z$. Nevertheless, such a form is still of interest, in particular because linear equations often occur in this form.
\begin{thm}\label{th:lag}
For the general linear equation \eqref{eq:2.3} in Laguerre-Forsyth form, the infinitesimal generators $X_n$ of $G_{\!S}$ and $X_n^0$ of $G_{\!c}$ have the following expressions, where $a_0, a_1, a_2,$ and $k_1$ are arbitrary constants, and $h$ an arbitrary function.
\begin{enumerate}
\item[{\rm a)}]
\begin{align}\label{eq:4.12}
\begin{split}
& X_n = (a_2 x^2 + a_1 x + a_0)\partial _x + \Big[ y \big[k_1 + \dfrac{n-1}{2}\big( 2a_2 x + a_1\big)\big] + h \Big] \partial _y \\
&\qquad \quad  +\sum _{k=0}^{n-3} \Big[  -(n-k)(2a_2x + a_1) c_k + a_2 (k+1)(k+1-n)c_{k+1}\\
&\qquad \quad +\delta_0^k \big[ -c_k \dfrac{h}{y} + \sum_{j=1}^{n-k}\binom{k+j}{j}\dfrac{h^{(j)}}{y}\big] \Big] \partial _{c_k}
\end{split}
\end{align}
\item[{\rm b)}]
\begin{align}\label{eq:4.13}
\begin{split}
& X_n^0 = \left(a_2 x^2 + a_1 x + a_0 \right)\partial _x\\
&\qquad \quad + \sum_{k=0}^{n-3}\left[-(n-k)(2a_2 x +a_1)c_k + a_2 (k+1)(k+1-n)c_{k+1} \right]\partial _{c_k}.
\end{split}
\end{align}
\end{enumerate}
\end{thm}
\begin{proof}
As in the proof of Theorem \ref{th:nor}, we only need to note that as the Laguerre-Forsyth form is a special case of the normal form, its generators $X_n$ and $X_n^0$ should be sought in the form \eqref{eq:4.8} and \eqref{eq:4.9}, respectively. More exactly, we only need to find the specific expression for the parameter $f$ of the equivalence transformation corresponding to the Laguerre-Forsyth form and substitute this into \eqref{eq:4.8} and \eqref{eq:4.9}, and to drop the term involving $c_{n-2}$ in the resulting expressions. It is well-known that the equivalence transformations of the Laguerre-Forsyth form of \eqref{eq:2.3} are invertible transformations of the form \eqref{eq:4.10} in which $T(z)$ is a linear fractional transformation. The corresponding infinitesimal generator is thus of the form \eqref{eq:4.10'}, in which $f(x)=a_2 x^2 + a_1 x + a_0$, for some arbitrary constants  $a_2, a_1$, and $a_0.$ This is the expression for $f$ which was to be found, and this completes the proof.
\end{proof}
Thanks to the ansatz \eqref{eq:4.4} a direct computation of $X_n$ and $X_n^0$ for equations of low orders up to seven has been performed and confirms the validity of the expressions given in the three preceding theorems. It should also be noted that unlike the case of equations in standard or in normal forms, the generator $X_n^0$ of $G_{\!c}$ in the case of the Laguerre-Forsyth form involves only a finite number of constant parameters. This means that the invariant functions for this form of the general linear equation are much easier to compute, as already noted by Forsyth \cite{forsy} who obtained an expression for them by a direct analysis.\par

    As noted earlier, for equations with a maximal symmetry algebra which are already expressed solely in terms of the two coefficients $c_{n-1}$ and $c_{n-2},$ to
obtain the corresponding infinitesimal generator $X^0,$ it suffices to substitute in the expression for $X_n^0$ corresponding to the general linear equation \eqref{eq:2.3} the corresponding characterizing equations which give an expression for the other coefficients solely in terms of $c_{n-1}$ and $c_{n-2}$ alone. For instance, for $n=4,$ the expression for $X_n^0$ corresponding to the normalized equation \eqref{eq:4.1} has, on account of \eqref{eq:2.14} and \eqref{eq:4.6}, an expression given by
\begin{equation}\label{eq:4.14}
\begin{split}
\xi  &= f \\
\phi_3^0 &= -c_3 f' - 4 g' + 6 f''\\
\phi_2^0 &= -2 c_2 f' - 3 c_3 g' + 3c_3 f'' - 6 g'' + 4 f^{(3)}\\
\phi_1^0 &= \frac{3}{8}f' (c_3^3 - 8 c_2' + 6 c_3 c_3' + 4 c_3'') - 3 c_3 g'' + c_3 f^{(3)} \\
         &\quad + c_2 \left(-\frac{3}{2} c_3 f' - 2 g' + f'' \right)- 4 g^{(3)} + f^{(4)}\\
\phi_0^0 &=  -\frac{1}{8} g' (8 c_2' - c_3 (-4 c_2 + c_3^2 + 6c_3')- 4 c_3'')- c_2 g'' -c_3 g^{(3)} \\
           &\quad  - g^{(4)} - \frac{1}{400} f' \Big[ 144 c_2^2 - 11c_3^4 - 288 c_3^2 c_3' - 8 c_2 (c_3^2 + 4 c_3') \\
           &\quad - 80 c_3 (5 c_2' - 7 c_3'') + 16 (21 c_3^{\prime 2} -30 c_2'' + 20c_3^{(3)})    \Big].
\end{split}
\end{equation}

\section{Concluding remarks}
We  reiterate the fact already mentioned that the symmetry properties obtained in this paper for linear equations also apply to the infinite dimensional class of nonlinear equations which are equivalent to a given linear equation admitting a maximal symmetry algebra. For instance, in the simplest case of the free fall equation $y''=0,$ an invertible point transformation of the form $x=f(z,w),\; y=g(z,w)$ shows that the most general class of second order (linear or nonlinear) equations admitting a maximal symmetry algebra has the form
\begin{equation*}
\begin{split}
 &f_z g_{z,z}-g_z f_{z,z}+w_z^3 \left(-g_w f_{w,w}+f_w g_{w,w}\right)\\
 &+w_z^2 \left(-g_z f_{w,w}-2 g_w f_{z,w}+f_z g_{w,w}+2 f_w g_{z,w}\right)\\
 &+w_z \left(-2 g_z f_{z,w}-g_w f_{z,z}+2 f_z g_{z,w}+f_w g_{z,z}\right)+\left(f_z g_w-f_w g_z\right) w_{z,z}=0.
\end{split}
\end{equation*}
Moreover, linearization methods under point transformations are available for {\sc ode}s of order up to three \cite{faz-leach, ibra-magri}, and this is very meaningful as for practical considerations most {\sc ode}s of physical relevance fall within this range.\par
One of the most interesting properties of linear equations with maximal symmetries is that their solution can be obtained by a very simple superposition formula from that of the second order source equation \cite{KM1}. More specifically, thanks to \eqref{eq:2.6}, any such equation can always be assumed to be in the normal reduced form \eqref{eq:2.4}. In particular, the corresponding second order source equation has the form $y'' + b y=0,$ for a certain function $b=b(x).$ If we let $u$ and $v$ be two linearly independent solutions of this source equation, then $n$ linearly independent solutions of an equation of the form \eqref{eq:2.4} with the same source equation are given by
$$y_k = u^k v^{n-1-k}, \qquad k=0,\dots,  n-1.$$
The latter fact can be used not only for finding analytic solutions of nonlinear equations, but also in the test of numerical schemes. Indeed, when testing a numerical scheme, it is always helpful to have an appropriate collection of nonlinear problems for which one or more explicit analytic solutions are available \cite{Bradie, Higham}.\par

The infinitesimal generators $X_n^0$ of the induced pseudo group of transformations $G_c$ found in Section \ref{sec:infigen} are of a more general interest. One of their main role is in the determination of the invariants (and semi-invariants) of the family of equations, and these functions can in turn be used for a complete classification of the given family of equations \cite{olv-fls2, mrz-clhyper}, thus reducing the study in each equivalence class to that of the canonical equation. For a much practical and immediate use, they are very efficient in testing whether a given function is an invariant of the related family of equation, and any given invariant of the family can also easily be used to test some necessary conditions of equivalence between two given equations.

\bibliographystyle{model1-num-names}

\end{document}